\documentclass{article}

\usepackage{
 amsmath,
 amsxtra,
 amsthm,
 amssymb,
 etex,
 mathrsfs,
 }
\usepackage[all]{xy}

\newtheorem{theorem}{Theorem}[section]
\newtheorem{lemma}[theorem]{Lemma}

\newtheorem{proposition}[theorem]{Proposition}
\newtheorem{corollary}[theorem]{Corollary}
\newtheorem*{thm2}{Theorem}
\newtheorem*{prop2}{Proposition}

\theoremstyle{definition}
\newtheorem{defn}[theorem]{Definition}
\newtheorem{remark}[theorem]{Remark}

\newcommand{\bd}{\begin{defn}}
\newcommand{\ed}{\end{defn}}
\newcommand{\bl}{\begin{lemma}}
\newcommand{\el}{\end{lemma}}
\newcommand{\bp}{\begin{proposition}}
\newcommand{\ep}{\end{proposition}}
\newcommand{\bt}{\begin{theorem}}
\newcommand{\et}{\end{theorem}}
\newcommand{\bc}{\begin{corollary}}
\newcommand{\ec}{\end{corollary}}
\newcommand{\br}{\begin{remark}}
\newcommand{\er}{\end{remark}}
\newcommand{\ba}{\begin{array}}
\newcommand{\ea}{\end{array}}
\newcommand{\bpf}{\begin{proof}}
\newcommand{\epf}{\end{proof}}

\newcommand{\Q}{\mathbb{Q}}
\newcommand{\Zp}{\mathbb{Z}_{p}}
\newcommand{\Qp}{\mathbb{Q}_{p}}
\newcommand{\Op}{\mathcal{O}}
\newcommand{\Ep}{E[p^{\infty}]}
\newcommand{\al}{\alpha}
\newcommand{\be}{\beta}
\newcommand{\Ga}{\Gamma}
\newcommand{\ga}{\gamma}
\newcommand{\e}{\varepsilon}
\newcommand{\La}{\Lambda}
\newcommand{\la}{\lambda}

\newcommand{\Si}{\Sigma}

\newcommand{\s}{\overrightarrow{s}}

\DeclareMathOperator{\Sel}{Sel} \DeclareMathOperator{\Gal}{Gal}
\DeclareMathOperator{\Hom}{Hom} \DeclareMathOperator{\rank}{rank}
\DeclareMathOperator{\corank}{corank}
 \DeclareMathOperator{\Tor}{Tor}

\newcommand{\cts}{\mathrm{cts}}

\newcommand{\M}{\mathfrak{M}}

\newcommand{\wE}{\widehat{E}}
\newcommand{\mK}{\mathcal{K}}
\newcommand{\mL}{\mathcal{L}}

\newcommand{\ot}{\otimes}
\newcommand{\ilim}{\displaystyle \mathop{\varinjlim}\limits}

\newcommand{\im}{\mathrm{im}\,}
\newcommand{\coker}{\mathrm{coker}\,}

\newcommand{\lra}{\longrightarrow}

\newcommand{\ps}[1]{[[ #1 ]]}

  \DeclareFontFamily{U}{wncy}{}
  \DeclareFontShape{U}{wncy}{m}{n}{<->wncyr10}{}
  \DeclareSymbolFont{mcy}{U}{wncy}{m}{n}
  \DeclareMathSymbol{\sha}{\mathord}{mcy}{"58}

\begin{document}
\title{On the cohomology of Kobayashi's plus/minus norm groups and applications}
 \author{
  Meng Fai Lim\footnote{School of Mathematics and Statistics $\&$ Hubei Key Laboratory of Mathematical Sciences,
Central China Normal University, Wuhan, 430079, P.R.China.
 E-mail: \texttt{limmf@mail.ccnu.edu.cn}} }
\date{}
\maketitle

\begin{abstract} \footnotesize
\noindent The plus and minus norm groups are constructed by Kobayashi as subgroups of the formal group of an elliptic curve with supersingular reduction, and they play an important role in Kobayashi's definition of the signed Selmer groups. In this paper, we study the cohomology of these plus and minus norm groups. In particular, we show that these plus and minus norm groups are cohomologically trivial. As an application of our analysis, we establish certain (quasi-)projectivity properties of the non-primitive mixed signed Selmer groups of an elliptic curve with good reduction at all primes above $p$. We then build on these projectivity results to derive a Kida formula for the signed Selmer groups under a slight weakening of the usual $\mu=0$ assumption, and study the integrality property of the characteristic element attached to the signed Selmer groups.


\end{abstract}

\section{Introduction}

Let $p$ be an odd prime. For now, let $E$ denote an elliptic curve over $\Qp$ with good supersingular reduction and $a_p=1 + p - |\tilde{E}(\mathbb{F}_p)| = 0$, where $\tilde{E}_v$ is the reduced curve of $E$. Let $\mathcal{K}$ be a finite unramified extension of $\Qp$, and let $\mK_\infty$ denote the cyclotomic $\Zp$-extension of $\mK$. Kobayashi \cite{Kob} defined the plus/minus norm groups $\wE^\pm(\mK_\infty)$, which are subgroups of the formal group of the said elliptic curve.
 These plus/minus norm groups are integral components in the definition of the signed Selmer groups. In particular, they come up in the local conditions at the supersingular primes in the definition of the Selmer groups, which were originally defined by Kobayashi \cite{Kob}, and later studied by other authors (for instances, see \cite{AL2,Kim07,KimPM,Kim14,KO, Leip-adic,LLMordell, LLAkashi, LZ}). In this paper, we will study the cohomology of the plus/minus norm groups, and our results are as follow.

\begin{prop2} [Proposition \ref{H invariant of plus-minus}]
For every subgroup $H$ of $\Gal(\mK/\Qp)$, we have
 \[ H^i\left(H, \widehat{E}^\pm(\mK_\infty)\otimes \Qp/\Zp\right) =\begin{cases} \widehat{E}^\pm(\mK_\infty^H)\otimes\Qp/\Zp,  & \mbox{if $i=0$}, \\
0, & \mbox{if $i \geq 1$}.\end{cases}\]
\end{prop2}

The $\Lambda$-structure of the $\widehat{E}^\pm(\mK_\infty)\otimes\Qp/\Zp$ has been well-studied (for instances, see \cite{Kim07, Kim14, KO, Kob}).
When $4$ divides $|\mK:\Qp|$, it has been observed by Kitajima-Otsuki \cite{KO} that $E^{+}(\mK_{\infty})\ot_{\Zp}\Qp/\Zp$ is not even cofree over $\La$, much less $\La[\Gal(\mK/\Qp]$. Therefore, Proposition \ref{H invariant of plus-minus} may seems rather surprising in view of these prior results. After this work is completed, Takenori Kataoka has pointed out to the author that the cohomological triviality of the plus/minus groups can be established by a slight refinement of the approach of Kitajima-Otsuki which applies also to the case when $4$ divides $|\mK:\Qp|$, and this was done in his thesis. (We thank him for making us aware of this.) We shall briefly say a bit on our proof which differs from the above approaches, and may be of independent interest. The proof requires a preliminary analysis of the formal group of $E$, where we have utilized a classical result of Noether \cite{No} on the Galois module structure of the ring of integers of unramified extensions and a deep result of Coates-Greenberg \cite{CG} on the formal groups of $E$ over infinitely ramified pro-$p$ extensions. Curiously, no prior knowledge of the explicit structural description of the plus/minus groups, as mentioned above, is required in the eventual proof.

We move to global consideration, which is the main motivation behind our study on the plus/minus norm groups. From now on, $E$ will denote an elliptic curve which is defined over a number field $F'$. Let $F$ be a finite extension of $F'$ and $K$ a finite Galois extension of $F$. The following assumptions will be in full force for our datum $(E, F/F', K)$ throughout.

\begin{itemize}
\item[(S1)] The elliptic curve $E$ has good reduction at all primes of $F'$ above $p$.

 \item[(S2)] There exists at least one prime of $F'$ above $p$ at which the elliptic curve $E$ has good supersingular reduction.

 \item[(S3)] For each prime $u$ of $F'$ above $p$ at which $E$ has supersingular reduction, we have

 (a) $F'_u=\Qp$ with $u$ being unramified in $K/F'$ (and hence in $F/F'$).

 (b) $a_u = 1 + p - |\tilde{E}_u(\mathbb{F}_p)| = 0$, where $\tilde{E}_u$ is the reduction of $E$ at $u$.
\end{itemize}

Let $\Si$ be a finite set of primes of $F$ which contains all the primes above $p$, all the ramified primes of $F/F'$ and $K/F$, the bad reduction primes of $E$ and the archimedean primes. Denote by $\Si_p$ (resp., $\Si_1$) the set of primes of $F$ above $p$ (resp. the set of primes of $F$ not dividing $p$). Write $\Si_p^{ord}$ (resp., $\Si_p^{ss}$) for the set of primes in $\Si_p$ at which $E$ has good ordinary reduction (resp., good supersingular reduction).

For each $v\in \Si_p^{ss}$, we can choose one of the two plus and minus norm subgroups. Fix such a collection of choices $\s =(s_v)_{v\in\Si^{ss}_p}\in \{\pm\}^{\Si^{ss}_p}$. By choosing the signs of each supersingular prime of $E$ over $K$ in a consistent manner with respect to $\s$ (see Subsection \ref{consistent sign}), we attach a signed Selmer group $\Sel^{\s}(E/K_{\infty})$ over the cyclotomic $\Zp$-extension $K_\infty$ of $K$. If $\Si_0$ is a subset of $\Si_1$, one can also define a non-primitive signed Selmer group $\Sel^{\s}_{\Si_0}(E/K_{\infty})$ (see Subsection \ref{non-primitive Sel subsection} for the precise definition). Every of these signed Selmer groups has a natural $\La[G]$-action, where
$\La=\Zp\ps{\Gal(K_{\infty}/K)}$ and $G=\Gal(K/F)$. We shall write $X^{\s}(E/K_{\infty})$ and $X^{\s}_{\Si_0}(E/K_{\infty})$ for the Pontryagin dual of $\Sel^{\s}(E/K_{\infty})$ and the Pontryagin dual of $\Sel^{\s}_{\Si_0}(E/K_{\infty})$ respectively.

Following Greenberg \cite{G11}, we set
\[\Phi_{K/F} = \{ v\in \Si_1~|~\mbox{the inertia degree of $v$ in $K/F$ is divisible by $p$} \}. \]
Our results on the non-primitive signed Selmer groups are as follow.

\begin{thm2}[Theorem \ref{main theorem1}]
Suppose that $(S1)-(S3)$ are satisfied and that $\Si_0$ contains $\Phi_{K/F}$. Assume that $X^{\s}(E/K_{\infty})$ is torsion over $\Zp\ps{\Ga}$ and that $\theta\Big(X^{\s}(E/K_{\infty})\Big)\leq 1$.
 Then $X^{\s}_{\Si_0}(E/K_{\infty})/X^{\s}_{\Si_0}(E/K_{\infty})[p]$ is quasi-projective as a $\Zp[G]$-module. Furthermore, if $X^{\s}(E/K_{\infty})$ is finitely generated over $\Zp$, then
  $X^{\s}_{\Si_0}(E/K_{\infty})$ is quasi-projective as a $\Zp[G]$-module.
\end{thm2}

Here $\theta(M)$ is the largest exponent of $p$ occurring in the structural decomposition of the $\Zp\ps{\Ga}$-module $M$.

\begin{thm2}[Theorem \ref{main theorem2}]
Suppose that $(S1)-(S3)$ are satisfied and that $\Si_0$ contains $\Phi_{K/F}$. Assume that  $\Sel^{\s}(E/K_{\infty})$ is cotorsion over $\Zp\ps{\Ga}$, and that for each $v\in \Si_p^{ord}$, either $v$ is non-anomalous for $E/K$ or $v$ is tamely ramified in $K/F$.
 Then $X^{\s}_{\Si_0}(E/K_{\infty})$ has a free resolution of length 1 as a $\La[G]$-module.
 Moreover, if one assumes further that $X^{\s}(E/K_{\infty})$ is finitely generated over $\Zp$, then $X^{\s}_{\Si_0}(E/K_{\infty})$ is projective as a $\Zp[G]$-module.
\end{thm2}

When the elliptic curve $E$ has good ordinary reduction at all primes above $p$, the above results were proved by Greenberg for the $p$-primary Selmer group in his monograph \cite{G11}. The proof of our theorem follows his strategy, where the essential difference lies in that we rely on our Proposition \ref{H invariant of plus-minus} to handle the local analysis at supersingular primes above $p$. Readers would have observe that we require the extension $K/F$ to be unramified at the supersingular primes (which is not required in the situation when $E$ has good ordinary reduction at all primes above $p$). This latter assumption is imposed on us by necessity due to that we do not have a nice enough decent theory for the plus/minus norm groups in ramified extension outside cyclotomic $\Zp$-extension (see \cite{LZ}).

The remainder of the paper is concerned with applications of the above theorems. As a start, we establish a Kida formula for the signed Selmer groups under a slight weakening of the usual $\mu=0$ assumption (see Proposition \ref{Kida formula}), which is inspired by an analogous result of Hachimori-Sharifi \cite{HSh} for the classical $\la$-invariants in cyclotomic $\Zp$-extensions of CM fields.
When the elliptic curve $E$ has good ordinary reduction at all primes above $p$, the Kida formula was first derived by Hachimori-Matsuno \cite{HM} building on an idea of Iwasawa \cite{Iw81} (also see \cite{G11}).  When the elliptic curve $E$ has good supersingular reduction at all primes above $p$, a Kida formula is also established for the signed Selmer groups in \cite{HL} and \cite{PW} by an argument via congruence of Galois representations.  Our approach differs from these and has the advantage of obtaining a Kida's formula under a weakening of the usual $\mu=0$ assumption, thus generalizing these results slightly. We should also mention that pertaining to this sort of theme, Alexandra Nichifor had worked out the connection between Kida's formula and projectivity for the classical Iwasawa modules in her thesis which unfortunately is unpublished. (The author likes to thank Bharathwaj Palvannan for pointing this out.)
Finally, we study certain integrality property of the characteristic elements of the non-primitive signed Selmer groups (see Proposition \ref{main theorem:char}). These characteristic elements are important objects of study and play a crucial role in the formulation of the main conjecture (see \cite{CFKSV, Leip-adic}). We should mention that for classical Iwasawa modules over totally real fields, the integrality property was recently studied by Nichifor-Palvannan in \cite{NP}. Our result can therefore be viewed as analogue of theirs.

We end the introductory section giving an outline of the paper. In
 Section \ref{projectivity properties}, we collect certain algebraic preliminaries that will be required in our subsequent discussion in the paper. In Section \ref{local calculations}, we analyze the cohomology of the formal groups and the plus/minus norm groups of a supersingular elliptic curve defined over $\Qp$. Section \ref{Selmer} is where we introduce the signed Selmer groups and the non-primitive signed Selmer groups. We shall also prove Theorems \ref{main theorem1} and \ref{main theorem2} here. Finally, in Section \ref{application section}, we discuss some applications of Theorems \ref{main theorem1} and \ref{main theorem2}.

\subsection*{Acknowledgement}
The author would like to thank Antonio Lei for many valuable discussions on the subject of the plus/minus norm groups. He would also like to thank Bharathwaj Palvannan for many insightful discussions and for answering many questions on his paper \cite{NP}. The author also liked to thank Takenori Kataoka and Andreas Nickel for their interest and comments on a previous version of the article. Finally, he thanks the referee for many helpful comments and suggestions on the paper. This research is supported by the National Natural Science Foundation of China under Grant No. 11550110172 and Grant No. 11771164.

\section{Projectivity properties} \label{projectivity properties}

In this section, we collect several results on the projectivity properties of modules, most of which can be found in the monograph \cite{G11}. For the convenience of the readers, we shall supply details to certain assertions that are essentially in \cite{G11} but perhaps not stated explicitly. In Subsection \ref{char subsection}, we also recall the notion of characteristic elements of Iwasawa modules following \cite{CFKSV}, and a certain maximal order following \cite{NP}, where the characteristic elements are expected to come from.
Throughout this section, $G$ will always denote a fixed finite group.

\subsection{$\Zp[G]$-modules}
As a start, we record the following criterion for a module to be projective over $\Zp[G]$.

\bp \label{proj1} Let $X$ be a $\Zp[G]$-module which is a free $\Zp$-module of finite rank. Write $S= \Hom(X,\Qp/\Zp)$. Suppose that $H^i(H,S)=0$ for $i=1,2$ and for every subgroup $H$ of $G$. Then $X$ is projective as a $\Zp[G]$-module.
\ep

\bpf See \cite[Proposition 2.1.1]{G11}.
\epf

As noted in \cite{G11}, for arithmetic applications, it is usually useful to work with the following weaker projectivity properties which we now recall.

\bd
A finitely generated $\Zp[G]$-module $X$ is said to be strictly quasi-projective if there exists a projective $\Zp[G]$-module $Y$ and a $\Zp[G]$-homomorphism $X\lra Y$ with finite kernel and cokernel.

The finitely generated $\Zp[G]$-module $X$ is said to be quasi-projective if there is a short exact sequence
\[ 0\lra X_1\lra X_2 \lra X\lra 0\]
of finitely generated $\Zp[G]$-modules with $X_1$ and $X_2$ being strictly quasi-projective as $\Zp[G]$-modules.

For alternative (and equivalent) definitions of strictly quasi-projective modules and  quasi-projective modules, we refer readers to \cite[Section 2]{G11}. \ed

For our purpose, it is useful to have a criterion for quasi-projectivity. To do this, we need to introduce certain notation and notions. Fix an integer $m\geq 1$ which is divisible by the order of all elements of $G$. Let $\mathcal{F}$ denote a finite extension of $\Qp$ which contains all $m$-th roots of unity, whose ring of integers is in turn denoted by $\Op$. For a cyclic subgroup $C$ of $G$, we write $C=PQ$, where $P$ is a $p$-group and $Q$ is a group with order prime to $p$.  Set $S=\Hom(X,\Qp/\Zp)$. For a $1$-dimensional character $\e$ of $Q$, we denote by $e_{\e}$ the idempotent for $\e$ in $\Op[Q]$. We shall then write $S^{\e} = e_{\e}(S\ot_{\Zp}\Op)$. In the subsequent discussion, we abbreviate $h_P(-)$ for the Herbrand quotient $|H^2(P,-)|/|H^1(P,-)|$.

\bp \label{proj2} Let $X$ be a finitely generated $\Zp[G]$-module.  Write $S= \Hom(X,\Qp/\Zp)$. Then $X$ is quasi-projective as a $\Zp[G]$-module if and only if for every cyclic subgroup $C=PQ$ of $G$ and every $1$-dimensional character character $\e$ of $Q$, we have $h_P(S^{\e})=1$.
\ep

\bpf See \cite[Proposition 2.1.3]{G11}.
\epf

The properties of a quasi-projectivity $\Zp[G]$-module are well-documented in \cite{G11}. For our purposes, we are mainly concerned with the following situation when $G$ is a finite $p$-group.

\bp \label{proj cor}
Suppose that $G$ is a finite $p$-group and that $X$ is a finitely generated quasi-projective $\Zp[G]$-module.
Then we have
\[ \rank_{\Zp}(X) = |G|\rank_{\Zp}(X_G).\]
\ep

This is essentially established in \cite[Sections 1 and 2]{G11} but not explicitly stated there. For the convenience of the reader, we give a proof here. To begin with, we have the following preparatory lemma.

\bl \label{proj lemma}
Suppose that $X$ is a finitely generated quasi-projective $\Zp[G]$-module.
Then for every cyclic normal subgroup $N$ of $G$, $X_N$ is quasi-projective as a $\Zp[G/N]$-module.
\el

\bpf
If $X$ is a projective $\Zp[G]$-module, then $X_N$ is clearly projective as a $\Zp[G/N]$-module. Now suppose that $X$ is strictly quasi-projective as a $\Zp[G]$-module. Hence we have a
$\Zp[G]$-homomorphism $X\lra Y$ with finite kernel and cokernel, and where $Y$ is some projective $\Zp[G]$-module. It is easy to see that the $\Zp[G]$-homomorphism $X\lra Y$ induces a $\Zp[G/N]$-homomorphism $X_N\lra Y_N$ with finite kernel and cokernel. But as noted at the beginning of the paragraph, $Y_N$ is projective as a $\Zp[G/N]$-module. Therefore, it follows that $X_N$ is strictly quasi-projective as a $\Zp[G/N]$-module.

Now suppose that $X$ is a finitely generated quasi-projective $\Zp[G]$-module. In other words, one has a short exact sequence
\[ 0\lra X_1\lra X_2 \lra X\lra 0\]
for some strictly quasi-projective $\Zp[G]$-modules $X_1$ and $X_2$. Taking $N$-invariant, we have the following exact sequence
\[ H_1(N,X)\lra (X_1)_N\stackrel{f}{\lra} (X_2)_N\lra X_N \lra 0.\]
Since $N$ is cyclic, we have that $H_1(N,X)$ is finite by Proposition \ref{proj2}. Combining this observation with that $(X_1)_N$ is a strictly quasi-projective $\Zp[G/N]$-module, it follows that $\im f$ is also strictly quasi-projective over $\Zp[G/N]$. Hence $X_N$ is a quasi-projective $\Zp[G/N]$-module.
\epf

We can now give the proof of Proposition \ref{proj cor}.

\bpf[Proof of Proposition \ref{proj cor}]
We shall first prove the proposition for a projective $\Zp[G]$-module. Since the ring $\Zp[G]$ is local, a projective $\Zp[G]$-module is necessarily free over $\Zp[G]$, and one clearly has the equality of the proposition in this situation.

Now, suppose that $X$ is strictly quasi-projective as a $\Zp[G]$-module. By definition, there exists some projective $\Zp[G]$-module $Y$ and a
$\Zp[G]$-homomorphism $X\lra Y$ with finite kernel and cokernel. The latter induces a $\Zp$-homomorphism $X_G\lra Y_G$ with finite kernel and cokernel. It then follows that
\[ \rank_{\Zp}(X) = \rank_{\Zp}(Y) = |G|\rank_{\Zp}(Y_G)= |G|\rank_{\Zp}(X_G),\]
where the middle equality follows from the discussion in the previous paragraph.

It remains to consider the situation that $X$ is a finitely generated quasi-projective $\Zp[G]$-module. For now, we assume that the group $G$ is cyclic of order $p$. Then one has a short exact sequence
\[ 0\lra X_1\lra X_2 \lra X\lra 0\]
for some strictly quasi-projective $\Zp[G]$-modules $X_1$ and $X_2$.
From this short exact sequence, we have
\[ H_1(G,X)\lra (X_1)_G\lra (X_2)_G\lra X_G \lra 0.\]
By the quasi-projective hypothesis and our assumption on $G$ being a cyclic group, we may apply Proposition \ref{proj2} to conclude that $H_1(G,X)$  is finite. Combining this observation with the above exact sequence, we obtain
\[ \rank_{\Zp}(X) = \rank_{\Zp}(X_2)-\rank_{\Zp}(X_1) = p\Big(\rank_{\Zp}(X_2)_G-\rank_{\Zp}(X_1)_G\Big)= p\rank_{\Zp}(X_G).\]
This proves the proposition for a quasi-projective $\Zp[G]$-module when $G$ is a cyclic group of order $p$.

For a general finite $p$-group $G$, let $N$ be a normal subgroup of $G$ which is cyclic of order $p$. Then by the previous paragraph, we have
\[ \rank_{\Zp}(X) = p\rank_{\Zp}(X_N).\]
By Lemma \ref{proj lemma}, $X_N$ is quasi-projective over $\Zp[G/N]$. Since $|G/N|<|G|$, it follows by induction that
\[ \rank_{\Zp}(X_N) = |G/N|\rank_{\Zp}(X_G).\]
Combining the two equalities, we obtain the conclusion of the proposition.
\epf

\subsection{$\lambda$-invariants}
We consider a slight refinement of Proposition \ref{proj2}, again following \cite{G11}. Let $\La$ denote the classical Iwasawa algebra $\Zp\ps{T}$. For a torsion $\La$-module $M$, there exist irreducible Weierstrass polynomials $f_j$, numbers $\al_i$, $\be_j$ and a homomorphism
\[ M \lra \La^r \oplus \Big(\bigoplus_{i=1}^s\La/p^{\al_i}\Big) \oplus \Big(\bigoplus_{j=1}^t\La/f_j^{\be_j}\Big) \]
with finite kernel and cokernel, where the numbers $\al_i$ $\be_j$ and the Weierstrass polynomials $f_j$ are determined by $M$ (cf. \cite[(5.3.8)]{NSW}). The $\mu$-invariant $\mu(M)$ (resp., $\la$-invariant $\la(M)$) of $M$ is defined to be $\sum_{i=1}^s\al_i$
(resp., $\sum_{j=1}^t\be_j\deg(f_j)$). We also set
\[\theta(M):=\max\{\al_i~|~i=1,..,s\}.\]
We now record two easy lemmas which can be deduced easily from the structural theorem and are left to the readers.

\bl \label{theta lemma} For a torsion $\La$-module, we have $\theta(M)\leq 1$ if and only if $M/M[p]$ is finitely generated over $\Zp$. Furthermore, in the event of such, we have
\[\la(M) = \rank_{\Zp}(M/M[p]).\]
\el

\bl \label{theta lemma2} Let $\phi: M\lra M'$ be a homomorphism of finitely generated torsion $\La$-modules, whose kernel and cokernel are finitely generated over $\Zp$. Then $\theta(M) = \theta(M')$.
\el

The next result is a slightly refined criterion for quasi-projectivity.

\bp \label{La proj2} Let $X$ be a finitely generated $\La[G]$-module which is torsion as a $\La$-module and satisfies $\theta(X)\leq 1$. Denote by $S$ the Pontryagin dual of $X$. Suppose that $H^1(N, S)$ and $H^2(N, S)$ are
finite for every cyclic $p$-subgroup $N$ of $G$. Then $X/X[p]$ is quasi-projective as a $\Zp[G]$-module if and only if for every cyclic subgroup $C=PQ$ of $G$ and every $1$-dimensional character character $\e$ of $Q$, we have $h_P(S^{\e})=1$.
\ep

\bpf See \cite[Proposition 2.2.1]{G11} or \cite[Section 2]{HSh}.
\epf

We end this subsection with the following refinement of Proposition \ref{proj cor}.

\bc \label{La proj cor}
Let $G$ be a finite $p$-group. Suppose that $X$ is a $\La[G]$-module which is torsion as a $\La$-module with $\theta(X)\leq 1$, and that $X/X[p]$ is quasi-projective as a $\Zp[G]$-module.
Then we have
\[ \la(X) = |G|\la(X_G).\]
\ec

\bpf
By Lemma \ref{theta lemma}, the module $X/X[p]$ is finitely generated over $\Zp$. Consequently,
we have
\[ \la(X) = \rank_{\Zp}(X/X[p]) = |G|\rank_{\Zp}\Big((X/X[p])_G\Big), \]
 where the first equality follows from Lemma \ref{theta lemma} and the second equality from Proposition \ref{proj cor}. On the other hand, we have an exact sequence
 \[X[p]_G\lra X_G \lra (X/X[p])_G \lra 0,\]
 where $X[p]_G$ is $p$-torsion. Thus, it follows that
 \[\la(X_G) =\rank_{\Zp}\Big((X/X[p])_G\Big),\]
 which yields the required conclusion of the corollary.
\epf

\subsection{Characteristic elements of $\La[G]$-modules} \label{char subsection}
As before, $G$ is a finite group, and $\La$ is the classical Iwasawa algebra $\Zp\ps{T}$. Following \cite{CFKSV, Kak11, Kak13}, we describe how to attach a characteristic element to a $\La[G]$-module which is torsion over $\La$.

Define the set
\[ S = \{ s\in \La[G]~|~ \La[G]/\La[G]s \mbox{ is finitely generated over } \Zp \}.\]
It follows from \cite[Theorem 2.4]{CFKSV} that $S$ is a left and right Ore set consisting of nonzero divisors in $\La[G]$.
Set $S^* = \cup_n p^nS$. By \cite[Proposition 2.3]{CFKSV}, a finitely generated $\La[G]$-module $M$ is annihilated by $S^*$ if and only if $M/M[p^{\infty}]$ is finitely generated over $\Zp$. Since our group $G$ is finite, the latter is equivalent to saying that $M$ is torsion over $\La$. In particular, the Grothendieck group of the category of finitely generated $\La[G]$-modules that are torsion over $\La$ identifies with the relative $K_0$ group $K_0(\La[G], \La[G]_{S^*})$. Henceforth, if $M$ is a finitely generated $\La[G]$-module which is torsion over $\La$, we may write $[M]$ for its class in $K_0(\La[G], \La[G]_{S^*})$.
On the other hand, by the proof of \cite[Lemma 2.1]{Kak11}, we have
\[ \La[G]_{S^*} \cong Q_\La[G],\]
where $Q_\La$ is the fraction field of $\La$. Under this identification, the localization sequence in $K$-theory yields the following exact sequence
\[ K_1(\La[G])\lra K_1( Q_\La[G]) \stackrel{\partial}{\lra} K_0(\La[G], Q_\La[G]).\]

\bl \label{surjective K}
The map $\partial$ in the above exact sequence is surjective.
\el

\bpf
When $G$ has no $p$-torsion, the surjectivity was established by Coates \textit{et al} \cite[Proposition 3.4]{CFKSV}. When $G$ has $p$-torsion, the surjectivity was established by Kakde in \cite[Lemma 1.5]{Kak11} for the Ore set $S$ rather than $S^*$. But it is not difficult to see that his proof goes through for $S^*$ with $G$ having $p$-torsion. Alternatively, one can also obtain the conclusion of lemma, appealing to the $K$-theory of Waldhausen categories as done by Witte \cite[Corollary 3.8]{Wi}.
\epf

In view of the preceding lemma, it makes sense to make the following definition.

\bd[Coates \textit{et al} \cite{CFKSV}] \label{Coates et al}
For a finitely generated $\La[G]$-module $M$ which is torsion over $\La$, a characteristic element for $M$ is an element $\xi_M\in K_1( Q_\La[G])$ such that $\partial(\xi_M) = [M]$.
\ed

\br
We emphasis that since our group $G$ is finite, we can attach characteristic element to $\La$-torsion $\La[G]$-module as above. In the case that $G$ is an infinite $p$-adic Lie group, one needs to impose a stronger hypothesis (the so-called $\M_H(G)$-conjecture) on the module $M$ to be able to define its characteristic element. We shall not need this in our paper but we do refer the interested readers to \cite{CFKSV, Kak11, Kak13, LZ} for discussion on this.
\er

To facilitate further discussion, we introduce a maximal $\La$-order in $Q_\La[G]$, following Nichifor-Palvannan \cite{NP}. By the Artin-Wedderburn theorem, there is an isomorphism
\[ \Qp[G]\cong \prod_iM_{m_i}(D_i)\]
of $\Qp$-algebras, where $M_{m_i}(D_i)$ is the ring of $m_i\times m_i$-matrices over $D_i$ and $D_i$ is a finite dimensional division algebra over $\Qp$. Denote by $\Op_{D_i}$ the maximal $\Zp$-order in $D_i$.
One then has the following commutative diagram
\[   \entrymodifiers={!! <0pt, .8ex>+} \SelectTips{eu}{}\xymatrix{
   \Zp[G]\ar[d] \ar[r]^{} & \Qp[G] \ar[d]^{\cong} \\
   \prod_iM_{m_i}(\Op_{D_i}) \ar[r]^{} &\prod_iM_{m_i}(D_i), } \]
   where the horizontal maps and leftmost vertical map are injective. This in turn induces the commutative diagram
   \[   \entrymodifiers={!! <0pt, .8ex>+} \SelectTips{eu}{}\xymatrix{
   \La[G]\ar[d] \ar[r]^{} & Q_\La[G] \ar[d]^{\cong} \\
   \prod_iM_{m_i}(\Op_{D_i}\ot_{\Zp}\La) \ar[r]^{} &\prod_iM_{m_i}(D_i\ot_{\Qp}Q_\La), } \]
 where the horizontal maps and leftmost vertical map are also injective.
 We shall set
 \[M_{\La[G]} := \prod_iM_{m_i}(\Op_{D_i}\ot_{\Zp}\La),\]
 which we view as a subring of $Q_\La[G]$ containing $\La[G]$. In fact, it follows from a combination of \cite[Proposition 2.11]{NP} and \cite[Theorem 10.5]{Re} that $M_{\La[G]}$ is a maximal $\La$-order of $Q_\La[G]$ which contains $\La[G]$.
We can now state the following result.

\bp \label{char proj 1}
Let $X$ be a finitely generated $\La[G]$-module which is $\La$-torsion. Suppose that $X$ has a free resolution of length one as a $\La[G]$-module. Then $\xi_X$ belongs to the image of the following natural map
\[M_{\La[G]}\cap Q_{\La}[G]^\times \lra K_1( Q_\La[G]).\]
\ep

\bpf
This was proved in \cite[Theorem 1]{NP} for a certain class of arithmetic Iwasawa modules. There they have an extra assumption on the validity of certain main conjecture \cite[Conjecture 3.2]{NP}. Going through the proof, one sees that the validity of the main conjecture is only used in attaching characteristic elements to their Iwasawa modules. However, for any module $X$ satisfying the hypothesis in our proposition, we have seen from Lemma \ref{surjective K} and Definition \ref{Coates et al} that one can always attach characteristic element to it. The remainder of the proof of \cite[Theorem 1]{NP} is purely algebraic relying on the free resolution hypothesis and an algebraic result \cite[Proposition 2.13]{NP}, and requires no additional input of the main conjecture. Going though the proof, one can see that this part of the argument carries over to prove the proposition.
\epf

As the preceding proposition requires a $\La[G]$-module to have a free resolution of length one, it will be useful to have the following criterion.

\bp \label{proj3} Let $X$ be a finitely generated $\La[G]$-module, where $\La=\Zp\ps{T}$. Suppose that $X$ is torsion as a $\La$-module and has no nonzero finite $\La$-submodule. Let $S = \Hom_{\cts}(X,\Qp/\Zp)$. Then $X$ has a free resolution of length one as a $\La[G]$-module if and only if $H^i(H,S)=0$ for $i=1,2$ and for every subgroup  $H$ of $G$.
\ep

\bpf See \cite[Proposition 2.4.1]{G11} or \cite[Proposition 4.1 and Remark 4.2]{NP}.
\epf

\section{Elliptic curve over $\Qp$ with supersingular reduction} \label{local calculations}

In this section, $E$ will denote an elliptic curve defined over $\Qp$ with supersingular reduction and $a_p=0$. We write $\wE$ for its formal group. Throughout, for an extension $L$ of $\Qp$, we write $\widehat{E}(L)$ for $\widehat{E}(\mathfrak{m}_L)$, where $\mathfrak{m}_L$ is the maximal ideal of the ring of integers of $L$. Let $\mK$ be a finite unramified extension of $\Qp$.
Denote by $\Q_{p,\infty}$ (resp., $\mK_{\infty}$) the cyclotomic $\Zp$-extension of $\Qp$ (resp., $\mK$). Since $\mK/\Qp$ is unramified and $\Q_{p,\infty}/\Qp$ is totally ramified, we have $\Q_{p,\infty}\cap \mK = \Qp$ and $\mK_{\infty} = \mK\Q_{p,\infty}$. Consequently, we have an isomorphism $\Gal(\mK_{\infty}/\Qp)\cong \Gal(\mK/\Qp)\times \Gal(\mK_{\infty}/\mK)$ of abelian groups.  Denote by $\mK_n$ the intermediate subextension of $\mK_{\infty}/\mK$ with $|\mK_n:\mK|=p^n$. The residue field of $\mK_n$ is in turn denoted by $k_n$.

\begin{lemma} \label{supersingular points}
 The formal group $\widehat{E}(\mK_n)$ has no $p$-torsion for every $n$. In particular, $E(\mK_n)$ has no $p$-torsion for every $n$.
 \end{lemma}

\bpf
The first assertion follows from \cite[Proposition 3.1]{KO} or \cite[Proposition 8.7]{Kob}. The second assertion follows from the first, the following short exact sequence
\[ 0\lra \widehat{E}(\mK_n)\lra E(\mK_n) \lra \widetilde{E}(k_n)\lra 0\]
and that $\widetilde{E}(k_n)$ has no $p$-torsion by the supersingular assumption. \epf

For subsequent discussion in this section, we shall write $G=\Gal(\mK/\Qp)$ which is also identified with $\Gal(\mK_\infty/\Q_{p,\infty})$. For a subgroup $H$ of $G$, we write $\mL:=\mL(H)$ for the fixed field of $\mK$ by $H$. Via the above identification, we may regard $H$ as a subgroup of $\Gal(\mK_{\infty}/\Q_{p,\infty})$. Under this setup, one can verify easily that the fixed field of $\mK_\infty$ by $H$ is precisely given by $\mL_\infty$, the cyclotomic $\Zp$-extension of $\mL$.

\subsection{Cohomology of formal groups} \label{coh of formal group subsec}
In this subsection, we will study the cohomology of the formal group. We begin with the following observation which builds on a classical result of Noether \cite{No}.

\bl \label{Noether free}
$\wE(\mK)$ is a free $\Zp[G]$-module of rank $1$. In particular, we have
\[H^i\big(H, \wE(\mK)\big) =0\]
for every subgroup $H$ of $G$ and $i\geq 1$.
\el

\bpf
Since $\mK$ is unramified over $\Qp$, and taking Lemma \ref{supersingular points} into account, the formal group logarithm gives an isomorphism $\wE(\mK)\cong p\Op_{\mK} \cong \Op_{\mK}$. Again, as $\mK$ is unramified over $\Qp$, we may apply a classical result of Noether \cite{No} to conclude that the latter is a free $\Zp[G]$-module of rank $1$. This proves the lemma.
\epf

\br
The cohomological triviality of $\wE(\mK)$ is a classical result (for instance, see \cite[Proposition 3.9]{CG} or \cite[Propositions 3.10 and 4.2]{EN}). In fact, this suffices for the application in this paper. But we have thought it interesting to note down that one can actually establish the $\Zp[G]$-freeness of $\wE(\mK)$ by appealing to Noether's result.
\er

\bp \label{H invariant of K}
For every subgroup $H$ of $G$, we have \[ H^i\left(H, \widehat{E}(\mK)\otimes \Qp/\Zp\right) =\begin{cases} \widehat{E}(\mL)\otimes\Qp/\Zp,  & \mbox{if $i=0$}, \\
0, & \mbox{if $i\geq 1$},\end{cases}\]
where $\mL = \mK^H$.
 \ep

\bpf
  By Lemma \ref{supersingular points}, we have $\Tor^{\Zp}_1(\widehat{E}(\mK),\Qp/\Zp)=0$ which in turn yields the following short exact sequence
   \[ 0\lra \widehat{E}(\mK)\lra \widehat{E}(\mK)\otimes\Qp \lra \widehat{E}(\mK) \otimes \Qp/\Zp\lra 0.\]
   In view of Lemma \ref{Noether free}, upon taking $H$-invariant, we obtain
\[ 0\lra \widehat{E}(\mL)\lra \widehat{E}(\mL)\otimes\Qp \lra \Big(\widehat{E}(\mK) \otimes \Qp/\Zp\Big)^H\lra 0\]
and
\[ H^l\Big(H,\widehat{E}(\mK)\otimes\Qp \Big)\cong H^l\Big(H,\widehat{E}(\mK)\otimes\Qp/\Zp \Big)\]
for $l\geq 1$. The isomorphism in the proposition for $i=0$ follows from the short exact sequence.  On the other hand, since $H$ is finite and $\widehat{E}(\mK)\otimes\Qp$ is torsionfree, we have $H^l\Big(H,\widehat{E}(\mK)\otimes\Qp \Big)=0$, or equivalently,
 $H^l\Big(H,\widehat{E}(\mK)\otimes\Qp/\Zp \Big) = 0$.
\epf

We now record the following result on the structure of $\wE(\mK_\infty)$ which is a partial analog of Lemma \ref{Noether free}.

\bl \label{H homology of formal}
We have $H^i\big(H, \wE(\mK_\infty)\big) =0$ for every subgroup $H$ of $G$ and
$i\geq 1$.
\el

\bpf
Recall that we write $\mL$ for the fixed field of $H$. Then as seen in the discussion before Subsection \ref{coh of formal group subsec}, the cyclotomic $\Zp$-extension  $\mL_{\infty}$ of $\mL$ is the fixed field of $\mK_\infty$ by $H$.
By a well-known result of Coates-Greenberg \cite[Corollary 3.2]{CG}, we have $H^i(\mK_\infty, \widehat{E}(\bar{\mK})) =0 = H^i(\mL_\infty, \widehat{E}(\bar{\mK}))$ for $i\geq 1$. Hence the spectral sequence
\[ H^i(H, H^j(\mK_\infty, \widehat{E}(\bar{\mK}))  \Longrightarrow H^{i+j}(\mL_\infty, \widehat{E}(\bar{\mK}))\] degenerates yielding the required conclusion of the lemma.
\epf

\bp \label{H invariant of Kcyc}
For every subgroup $H$ of $G$, we have \[ H^i\left(H, \widehat{E}(\mK_\infty)\otimes \Qp/\Zp\right) =\begin{cases} \widehat{E}(\mL_\infty)\otimes\Qp/\Zp,  & \mbox{if $i=0$}, \\
0, & \mbox{if $i\geq 1$},\end{cases}\]
where $\mL_\infty = \mK_\infty^H$.
 \ep

\bpf
  This has a similar proof to that in Proposition \ref{H invariant of K}, where we make use of Lemma \ref{H homology of formal} in place of Lemma \ref{Noether free}.
\epf

\subsection{Cohomology of plus/minus norm groups}

In this subsection, we study the cohomology of the plus/minus norm groups.
Following \cite{KimPM,Kim14,KO,Kob}, we define
the plus and minus norm groups
\[\widehat{E}^+(\mK_n) = \left\{ P\in \widehat{E}(\mK_n)~:~\mathrm{tr}_{n/m+1}(P)\in \wE(\mK_l), 2\mid m, 0\leq m\leq n-1\right\}, \]
\[\widehat{E}^-(\mK_n) = \left\{ P\in \widehat{E}(\mK_n)~:~\mathrm{tr}_{n/m+1}(P)\in \wE(\mK_l), 2\nmid m, 0\leq m \leq n-1\right\}, \]
where $\mathrm{tr}_{n/m+1}: \widehat{E}(\mK_n) \lra \widehat{E}(\mK_{m+1})$ denotes the trace map. These groups are related to the formal groups in the following way.

\bl \label{plus/minus and formal}
There is a short exact sequence \[ 0 \lra  \widehat{E}(\mK) \lra \widehat{E}^+(\mK_\infty)\oplus \widehat{E}^-(\mK_\infty) \lra
\widehat{E}(\mK_\infty) \lra 0.\]
\el

\bpf
See \cite[Proposition 2.6]{Kim14} or \cite[Proposition 8.12]{Kob}.
\epf

We can now prove the following proposition as stated in the introductory section.

\bp \label{H invariant of plus-minus}
For every subgroup $H$ of $G$, we have
 \[ H^i\left(H, \widehat{E}^\pm(\mK_\infty)\otimes \Qp/\Zp\right) =\begin{cases} \widehat{E}^\pm(\mL_{\infty})\otimes\Qp/\Zp,  & \mbox{if $i=0$}, \\
0, & \mbox{if $i \geq 1$}.\end{cases}\]
\ep

\bpf
 In view of Lemmas \ref{supersingular points} and \ref{plus/minus and formal}, we have the following short exact sequences
\[0 \lra  \widehat{E}(\mL)\otimes\Qp/\Zp \lra \big(\widehat{E}^+(\mL_\infty)\otimes\Qp/\Zp\big)\oplus \big( \widehat{E}^-(\mL_\infty)\otimes\Qp/\Zp \big)\lra
\widehat{E}(\mL_\infty)\otimes\Qp/\Zp \lra 0,\]
\[0 \lra  \widehat{E}(\mK)\otimes\Qp/\Zp \lra \big(\widehat{E}^+(\mK_\infty)\otimes\Qp/\Zp\big)\oplus \big( \widehat{E}^-(\mK_\infty)\otimes\Qp/\Zp \big)\lra
\widehat{E}(\mK_\infty)\otimes\Qp/\Zp \lra 0,\]
which in turn fit into the following diagram
\[   \entrymodifiers={!! <0pt, .8ex>+} \SelectTips{eu}{}\xymatrix{
    0 \ar[r]^{} & \widehat{E}(\mL)\otimes\Qp/\Zp \ar[d] \ar[r] &  \big(\widehat{E}^+(\mL_\infty)\otimes\Qp/\Zp\big)\oplus \big( \widehat{E}^-(\mL_\infty)\otimes\Qp/\Zp \big)
    \ar[d] \ar[r] & \widehat{E}(\mL_\infty)\otimes\Qp/\Zp\ar[d]^{}  \ar[r] & 0\\
    0 \ar[r]^{} & \big(\widehat{E}(\mK)\otimes\Qp/\Zp\big)^H \ar[r]^{} & \big(\widehat{E}^+(\mK_\infty)\otimes\Qp/\Zp\big)^H\oplus \big( \widehat{E}^-(\mK_\infty)\otimes\Qp/\Zp \big)^H \ar[r] & \big(\widehat{E}(\mK_\infty)\otimes\Qp/\Zp\big)^H\ar[r]
    & \cdots } \]
Since the leftmost and rightmost vertical maps are isomorphisms by Propositions \ref{H invariant of K} and \ref{H invariant of Kcyc}, so is the middle map by a snake lemma argument. This yields the isomorphism of the proposition for $i=0$. Now the bottom sequence of the diagram continues in the form
\[ H^i(H, \widehat{E}(\mK)\otimes\Qp/\Zp) \lra H^i(H,\widehat{E}^+(\mK_\infty)\otimes\Qp/\Zp)\oplus H^i(H,\widehat{E}^-(\mK_\infty)\otimes\Qp/\Zp) \lra H^i(H,\widehat{E}(\mK_\infty)\otimes\Qp/\Zp).\]
Again, by virtue of Propositions \ref{H invariant of K} and \ref{H invariant of Kcyc},
we obtain the desired vanishing for $i\geq 1$.
\epf

Recall that $\widehat{E}^\pm(\mK_n)\otimes\Qp/\Zp$ injects into $H^1(\mK_n,\Ep)$ via the Kummer map (cf. \cite[Lemma 8.17]{Kob}). Thus, it makes sense to speak of $\displaystyle \frac{H^1(\mK_\infty, \Ep)}{\widehat{E}^\pm(\mK_\infty)\otimes\Qp/\Zp}$. We may now state and prove the following.

\bp \label{H-cohomology of supersingular}
For every subgroup $H$ of $G$, we have
\[ H^i\left(H, \frac{H^1(\mK_\infty, \Ep)}{\widehat{E}^\pm(\mK_\infty)\otimes\Qp/\Zp}\right) =\begin{cases} \displaystyle\frac{H^1(\mL_\infty, \Ep)}{\widehat{E}^\pm(\mL_\infty)\otimes\Qp/\Zp},  & \mbox{if $i=0$}, \\
0, & \mbox{if $i \geq 1$}.\end{cases}\]
\ep

\begin{proof}
 Consider the spectral sequence
  \[ H^i\big(H, H^j(\mK_\infty, \Ep)\big)\Longrightarrow H^{i+j}(\mL_\infty,\Ep). \]
  By \cite[Theorem 7.1.8(i)]{NSW}, $H^l(\mK_\infty, \Ep) = 0 =H^l(\mL_\infty, \Ep)$ for $l\geq 2$. Also, recall that $E(\mK_\infty)[p^\infty] =0$ by Lemma \ref{supersingular points}. Hence the spectral sequence degenerates to yield
  \[ H^i\left(H, H^1(\mK_\infty, \Ep)\right) =\begin{cases} H^1(\mL_\infty, \Ep),  & \mbox{if $i=0$}, \\
0, & \mbox{if $i\geq 1$}.\end{cases}\]
The conclusion of the corollary now follows from combining this latter observation with an analysis of the $H$-cohomology exact sequence of
\[ 0\lra \widehat{E}^\pm(\mK_\infty)\otimes\Qp/\Zp \lra H^1(\mK_\infty, \Ep)\lra \frac{H^1(\mK_\infty, \Ep)}{\widehat{E}^\pm(\mK_\infty)\otimes\Qp/\Zp}\lra 0\]
and taking Proposition \ref{H invariant of plus-minus} into account.
\end{proof}

It follows immediately from Proposition \ref{H-cohomology of supersingular} that we have the following.

\bc \label{H-cohomology of supersingular2}
For every cyclic subgroup $C=PQ$ of $G$ and every $1$-dimensional character $\e$ of $Q$, we have
\[ h_P\left(\Big(\frac{H^1(\mK_\infty, \Ep)}{\widehat{E}^\pm(\mK_\infty)\otimes\Qp/\Zp}\Big)^{\e} \right)=1.\]
\ec

\bpf
Write $M=\frac{H^1(\mK_\infty, \Ep)}{\widehat{E}^\pm(\mK_\infty)\otimes\Qp/\Zp}$. By Proposition \ref{H-cohomology of supersingular}, we have $H^i(P, M) = 0$ for $i\geq 1$. Let $\mathcal{F}$ denote the ring of integers of a finite extension of $\Qp$ which contains all $m$-th roots of unity, where $m$ is divisible by the order of all elements of $G$. We clearly still have $H^i(P, M\ot\Op) = 0$ for $i\geq 1$. On the other hand, since $M\ot\Op = \oplus_{\e}e_{\e}(M\ot_{\Zp}\Op)$, it follows that we have $H^i(P, e_{\e}(M\ot_{\Zp}\Op)) = 0$ for $i\geq 1$. The conclusion of the corollary is now immediate from this.
\epf

\br \label{relative +-}
One can actually have relative versions of the results in this section in the following sense. Let $\mK$ and $\mathcal{M}$ be two finite unramified extensions of $\Qp$ with $\mK\subseteq \mathcal{M}$.
Then every subgroup $H'$ of $\Gal(\mathcal{M}/\mK)$, the same argument in this section can prove the following version of Propositions \ref{H invariant of plus-minus} and \ref{H-cohomology of supersingular}:
 \[ H^i\left(H', \widehat{E}^\pm(\mathcal{M}_\infty)\otimes \Qp/\Zp\right) =\begin{cases} \widehat{E}^\pm(\mathcal{M}_\infty^{H'})\otimes\Qp/\Zp,  & \mbox{if $i=0$}, \\
0, & \mbox{if $i \geq 1$},\end{cases}\]
and
\[ H^i\left(H', \frac{H^1(\mathcal{M}_\infty, \Ep)}{\widehat{E}^\pm(\mathcal{M}_\infty)\otimes\Qp/\Zp}\right) =\begin{cases} \displaystyle\frac{H^1(\mathcal{M}_\infty^{H'}, \Ep)}{\widehat{E}^\pm(\mathcal{M}_\infty^{H'})\otimes\Qp/\Zp},  & \mbox{if $i=0$}, \\
0, & \mbox{if $i \geq 1$}.\end{cases}\]
Similarly, one can also state relative versions of Proposition \ref{H invariant of K}, Proposition \ref{H invariant of Kcyc} and Corollary \ref{H-cohomology of supersingular2} which we leave to the readers to fill in.
\er

\section{Generalities on signed Selmer groups} \label{Selmer}

We now turn to the global situation, where we begin by fixing some notation and standing assumptions that will be adhered throughout
the section. Let $E$ be an elliptic curve defined over a number field $F'$ and $F$ a finite extension of $F'$. Fix a finite Galois extension $K$ of $F$. The following assumptions will be in full force for our data $(E, F/F', K)$ throughout.

\begin{itemize}
\item[(S1)] The elliptic curve $E$ has good reduction at all primes of $F'$ above $p$.

 \item[(S2)] There exists at least one prime of $F'$ above $p$ at which the elliptic curve $E$ has good supersingular reduction.

 \item[(S3)] For each $u$ of $F'$ above $p$ at which $E$ has supersingular reduction, we have

 (a) $F'_u=\Qp$ and $u$ is unramified in $K/F'$.

 (b) $a_u = 1 + p - |\tilde{E}_u(\mathbb{F}_p)| = 0$, where $\tilde{E}_u$ is the reduction of $E$ at $u$.
\end{itemize}

Let $\Si$ be a finite set of primes of $F$ which contains all the primes above $p$, all the ramified primes of $F/F'$ and $K/F$, the bad reduction primes of $E$ and the archimedean primes. Denote by $\Si_p$ the set of primes of $K$ above $p$. Write $\Si_p^{ord}$ (resp., $\Si_p^{ss}$) for the set of primes in $\Si_p$ at which $E$ has good ordinary reduction (resp., good supersingular reduction). We shall also write $\Si_1= \Si-\Si_p$ for the set of primes of $F$ which are not divisible by $p$.
For any subset $R$ of $\Si$ and any extension $\mathcal{F}$ of $F$, we shall write $R(\mathcal{F})$ for the set of primes of $\mathcal{F}$ above $R$.

Denote by $K_\infty$ the cyclotomic $\Zp$-extension of $K$ and $K_n$ the intermediate subfield of $K_{\infty}/K$ with $|K_n:K|=p^n$. Let $L$ be an intermediate subextension of $K/F$. By (S2)-(S3), $K\cap L_\infty= L$, where $L_\infty$ is the cyclotomic $\Zp$-extension of $L$. Furthermore, denoting by $L_n$ the intermediate intermediate subfield of $L_{\infty}/L$ with $|L_n:L|=p^n$, we have $K_n\cap L_\infty = L_n$ and $K_n = L_nK$. Thus, we have $\Gal(K_{\infty}/L)\cong \Gal(L_{\infty}/L)\times \Gal(K/L)$. We write $\Ga$ for $\Gal(F_\infty/F)$ which is identified with $\Gal(L_\infty/L)$ and, in particular, $\Gal(K_\infty/K)$. We shall also write $H_L =\Gal(K/L)$ which is identified with $\Gal(K_\infty/L_\infty)$. Finally, set $G=H_F = \Gal(K/F)$.

\subsection{A consistent choice of signs} \label{consistent sign}

We now define our signed Selmer group. As we need to compare signed Selmer groups over different extensions, our signs for the Selmer group have to be chosen in a consistent way which we now describe. We shall always fix $\overrightarrow{s}=(s_v)_{v\in \Si_p^{ss}}\in\{\pm\}^{\Si_p^{ss}}$.
Let $L$ be an intermediate subextension of $F/K$. By (S2)-(S3), for each prime $u$ in $\Si_p^{ss}(L)$, there is a unique prime of $L_n$ lying above the said prime which, by abuse of notation, is still denoted by $u$. For each $u\in \Si_p^{ss}(L_n)$, we set $s_u=s_v$, where $v$ is a prime of $F$ below $u$.
The $\overrightarrow{s}$-signed Selmer groups $\Sel^{\s}(E/L_n)$ over $L_n$ are then defined by
\[ \ker \Bigg(H^1(G_\Si(L_n),\Ep)\lra \bigoplus_{u\in \Si_p^{ord}(L_n)}\frac{H^1(L_{n,u},\Ep)}{E(L_{n,u})\ot\Qp/\Zp}\times\bigoplus_{u\in \Si_p^{ss}(L_n)}\frac{H^1(L_{n,v},\Ep)}{E^{s_u}(L_{n,u})\ot\Qp/\Zp}\]
\[\times\bigoplus_{v\in \Si_1(L_n)}H^1(L_{n,u},\Ep) \Bigg).\]
We set $\Sel^{\s}(E/L_{\infty})=\ilim_n \Sel^{\overrightarrow{s}}(E/L_n)$. In particular, we have
\[ \Sel^{\s}(E/L_{\infty})= \ker \Bigg(H^1(G_\Si(L_\infty),\Ep)\lra \bigoplus_{u\in \Si_p^{ord}(L_\infty)}\frac{H^1(L_{\infty,u},\Ep)}{E(L_{\infty,u})\ot\Qp/\Zp}\times\bigoplus_{u\in \Si_p^{ss}(L_\infty)}\frac{H^1(L_{\infty,v},\Ep)}{E^{s_u}(L_{\infty,u})\ot\Qp/\Zp}\]
\[\times\bigoplus_{v\in \Si_1(L_\infty)}H^1(L_{\infty,u},\Ep) \Bigg).\]
To simplify notation, for each $u\in \Si(L_{\infty})$, we shall write $J_u(E/L_{\infty})$ for $\displaystyle\frac{H^1(L_{\infty,u},\Ep)}{E(L_{\infty,u})\ot\Qp/\Zp}$, $\displaystyle\frac{H^1(L_{\infty,u},\Ep)}{E^{s_v}(L_{\infty,u})\ot\Qp/\Zp}$ or $H^1(L_{\infty,u},\Ep)$ accordingly to $u$ belonging to $\Si_p^{ord}(L_{\infty})$, $\Si_p^{ss}(L_{\infty})$  or $\Si_1(L_{\infty})$.

By the manner of the choice of signs, we then have a natural map
\[ \Sel^{\s}(E/L_{\infty}) \lra \Sel^{\s}(E/K_{\infty})^{H_L}\]
induced by the restriction map on cohomology.

\bl \label{descent lemma}
The map
\[ \Sel^{\s}(E/L_{\infty}) \lra \Sel^{\s}(E/K_{\infty})^{H_L}\]
is injective with finite cokernel. \el

\bpf
We have the following commutative diagram
\[   \entrymodifiers={!! <0pt, .8ex>+} \SelectTips{eu}{}\xymatrix{
    0 \ar[r]^{} &\Sel^{\s}(E/L_{\infty})  \ar[d] \ar[r] &  H^1\big(G_{\Sigma}(L_{\infty}),\Ep\big)
    \ar[d]^{h_L} \ar[r] & \displaystyle \bigoplus_{u\in\Sigma(L_{\infty})}J_u(E/L_{\infty})\ar[d]^{l=\oplus l_u}  \\
    0 \ar[r]^{} & \Sel^{\s}(E/K_{\infty})^{H_L} \ar[r]^{} & H^1\big(G_{\Sigma}(K_\infty),\Ep\big)^{H_L} \ar[r] &
    \left(\displaystyle\bigoplus_{w\in\Si(K_\infty)}J_w(E/F_{\infty})\right)^{H_L}  } \]
    with exact rows. For $w\in \Si^{ss}_p(K_\infty)$, we have $E(K_{\infty,w})[p^{\infty}]=0$ by Lemma \ref{supersingular points}, which in turn implies that $E(F_{\infty})[p^{\infty}]=0$ noting (S2).
Combining this observation with a Hochschild-Serre spectral sequence argument, we see that $h_L$ is an isomorphism. For $u\in \Si^{ss}_p(L_{\infty})$, it follows from Proposition \ref{H-cohomology of supersingular} and Remark \ref{relative +-} that $l_u$ is an isomorphism. Finally, it follows from the proof of \cite[Lemma 3.3]{HM} that $\ker l_u$ is finite for the remaining primes. The assertion of the lemma is now an immediate consequence of these observations.
\epf

\medskip \noindent \textbf{Conjecture.} Write $X^{\s}(E/K_{\infty})$ for the Pontryagin dual of $\Sel^{\overrightarrow{s}}(E/K_{\infty})$. Then $X^{\s}(E/K_{\infty})$ is a torsion $\Zp\ps{\Ga}$-module, where $\Ga=\Gal(K_{\infty}/K)\cong\Zp$.

\medskip
When $E$ has good ordinary reduction at all primes above $p$, the above conjecture is precisely Mazur's conjecture \cite{Maz} which is known in the case when $E$ is defined over $\Q$ and $F$ an abelian extension of $\Q$ (see \cite{K}). For an elliptic curve over $\Q$ with good supersingular reduction at $p$, this conjecture was established by Kobayashi (cf.\ \cite{Kob}; also see \cite{BL} for some recent progress on this conjecture). We now record an equivalent characterization of this property.

\bp \label{torsion surjective H2}
The module $X^{\s}(E/K_{\infty})$ is $\Zp\ps{\Ga}$-torsion if and only if we have $H^2(G_S(K_{\infty}),\Ep)=0$ and the following short exact sequence
\[0\lra \Sel^{\overrightarrow{s}}(E/K_{\infty}) \lra H^1(G_\Si(K_\infty),\Ep)\lra \bigoplus_{w\in \Si(K_\infty)}J_w(E/K_{\infty}) \lra 0. \]
\ep

\bpf
See \cite[Proposition 3.3.1]{AL2} or \cite[Proposition 2.7]{LLMordell}.
\epf

\br  \label{tor remark}
If $\Sel^{\overrightarrow{s}}(E/K_{\infty})$ is a cotorsion $\Zp\ps{\Ga}$-module, then so is $\Sel^{\overrightarrow{s}}(E/L_{\infty})$ for every intermediate subextension $L$ of $K/F$ by Lemma \ref{descent lemma}. Hence Proposition \ref{torsion surjective H2} also applies to $\Sel^{\overrightarrow{s}}(E/L_{\infty})$.
\er

We give two more consequences of $X^{\overrightarrow{s}}(E/K_{\infty})$ being torsion $\Zp\ps{\Ga}$-module.

\bp \label{tor cohomology}
Suppose that $\Sel^{\overrightarrow{s}}(E/K_{\infty})$ is a cotorsion $\Zp\ps{\Ga}$-module. Then  the cohomology group
$H^i(H, \Sel^{\s}(E/K_{\infty}))$ is finite for every subgroup $H$ of $G$ and $i\geq 1$.
\ep

\bpf
As before, we write $L$ for the fixed field of $H$. In view of Proposition \ref{torsion surjective H2}, we have the following commutative diagram
\[   \entrymodifiers={!! <0pt, .8ex>+} \SelectTips{eu}{}\xymatrix{
    0 \ar[r]^{} &\Sel^{\s}(E/L_{\infty})  \ar[d] \ar[r] &  H^1\big(G_{\Sigma}(L_{\infty}),\Ep\big)
    \ar[d]^{h_L} \ar[r] & \displaystyle \bigoplus_{u\in\Sigma(L_{\infty})}J_u(E/L_{\infty})\ar[d]^{l=\oplus l_u}\ar[r]^{}  &0 \\
    0 \ar[r]^{} & \Sel^{\s}(E/K_{\infty})^{H} \ar[r]^{} & H^1\big(G_{\Sigma}(K_\infty),\Ep\big)^{H} \ar[r] &
    \left(\displaystyle \bigoplus_{w\in\Sigma(K_{\infty})}J_w(E/K_{\infty})\right)^{H} \ar[r]^{} & \cdots } \]
    with exact rows, and the bottom row continues as

\[ H^1\Big(H, \Sel^{\s}(E/K_{\infty})\Big) \lra H^1\Big(H, H^1(G_{\Si}(K_\infty), \Ep)\Big) \lra H^1\left(H,\displaystyle \bigoplus_{w\in\Sigma(K_{\infty})}J_w(E/K_{\infty})\right)\lra \cdots \]

Now consider the spectral sequence
\[ H^i(H_L, H^j(G_{\Si}(K_{\infty}),\Ep))\Longrightarrow H^{i+j}(G_{\Si}(L_{\infty}),\Ep).\]
As seen in the proof of Lemma \ref{descent lemma}, $H^0(G_{\Si}(K_{\infty}),\Ep)=0 =H^0(G_{\Si}(L_{\infty}),\Ep)$. By Proposition \ref{torsion surjective H2} and Remark \ref{tor remark}, we have $H^2(G_{\Si}(K_{\infty}),\Ep)=0 =H^2(G_{\Si}(L_{\infty}),\Ep)$. Hence the spectral sequence degenerates yielding
\[ H^i(H_L, H^1(G_{\Si}(K_{\infty}),\Ep)) =\begin{cases} H^1(G_{\Si}(L_{\infty}),\Ep),  & \mbox{if $i=0$}, \\
0, & \mbox{if $i\geq 1$}.\end{cases} \]
As seen in the proof of Lemma \ref{descent lemma},
$h_L$ is an isomorphism. For $w\in \Si_p^{ss}(K_\infty)$, it follows from Proposition \ref{H-cohomology of supersingular} and Remark \ref{relative +-} that $\coker l_u=0$ and $H^i\Big(H, J_w(E/K_{\infty})\Big) =0$. For the remaining primes, it follows from the discussion in \cite[Section 4]{HM} that $\coker l_u$ and $H^i\Big(H, J_w(E/K_{\infty})\Big)$ are finite. Putting all these observations into the long exact sequence of the bottom row, we obtain the conclusion of the corollary.
\epf

\bp \label{no finite submodules}
Suppose that $X^{\overrightarrow{s}}(E/K_{\infty})$ is a torsion $\Zp\ps{\Ga}$-module. Then $X^{\overrightarrow{s}}(E/K_{\infty})$ has no nonzero finite $\Zp\ps{\Ga}$-submodules.
\ep

\bpf
 We follow the approach of \cite[Proposition 4.14]{G99} and \cite[Theorem 3.14]{KimPM}. Fix an isomorphism $\kappa:\Ga\cong 1+p\Zp$. For each $s$, we write $\Zp(s)$ for the abelian group $\Zp$ with a $\Ga$-action given by $\ga\cdot x = \kappa(\ga)x$. For every $\Ga$-module $M$, write $M(s) = M\ot \Zp(s)$, where $\Ga$ acts diagonally. As in \cite[Proposition 4.14]{G99} or \cite[Theorem 3.14]{KimPM}, one can choose an $s$ such that the twisted Selmer group $\Sel(E/K)(s)$ is finite. Then we have the following commutative diagram
 \[   \entrymodifiers={!! <0pt, .8ex>+} \SelectTips{eu}{}\xymatrix{
    0 \ar[r]^{} &\Sel(E/K)(s)  \ar[d] \ar[r] &  H^1\big(G_{\Sigma}(K),\Ep(s)\big)
    \ar[d]^{h'} \ar[r] & \displaystyle \bigoplus_{u\in\Sigma(K)}J_u(E/K)(s)\ar[d]^{\oplus g_u}\ar[r]^{}  &0 \\
    0 \ar[r]^{} & \Big(\Sel^{\s}(E/K_{\infty})(s)\Big)^{\Ga} \ar[r]^{} & \Big(H^1\big(G_{\Sigma}(K_\infty),\Ep\big)(s)\Big)^{\Ga} \ar[r] &
    \left(\displaystyle \bigoplus_{w\in\Sigma(K_{\infty})}J_w(E/K_{\infty})(s)\right)^{\Ga} \ar[r]^{} & \cdots } \]
    with exact rows. Here the surjectivity of the rightmost map in the top row follows from the finiteness of $\Sel(E/K)(s)$ by a similar argument to that in \cite[Proposition 3.8]{KimPM}. By a similar argument to that in Proposition \ref{tor cohomology}, the map $h'$ is bijective. For $u\in\Si^{ss}_p$, the map $g_u$ is given by
     \[g_u : \frac{H^1(K_u, \Ep)}{E(K_u)\ot \Qp/\Zp} \lra \left(\frac{H^1(K_{\infty,w}, \Ep)}{E^{\pm}(K_{\infty,w})\ot \Qp/\Zp}\right)^{\Ga}. \]
    It follows from \cite[Proposition 3.32]{KO} that the term $\left(\frac{H^1(K_{\infty,w}, \Ep)}{E^{\pm}(K_{\infty,w})\ot \Qp/\Zp}\right)^{\Ga}$ is a cofree $\Zp$-module and hence $\coker g_u$ is cofree over $\Zp$.
      (Note:\ the map $g_u$ need not be surjective; see \cite[Section 3]{AL}.) For each of the remaining primes, the map $g_u$ is surjective, since $\Ga$ has $p$-cohomological dimension one (for instance, see \cite[Proposition 4.14]{G99}). Hence it follows that we have an exact sequence
    \[ 0\lra \bigoplus_{u\in\Si^{ss}_p(K)}\coker g_u \lra H^1\Big(\Ga, \Sel^{\s}(E/K_{\infty})(s)\Big) \lra
    H^1\Big(\Ga, H^1(G_\Si(K_{\infty}),\Ep)(s)\Big).\]
    Via a similar argument to that in \cite[Lemma 2.5]{AL}, we see that $H^1\Big(\Ga, H^1(G_\Si(K_{\infty}),\Ep)(s)\Big) =0$ using the fact that $\Sel(E/K)(s)$ is finite. Now, it follows from the above exact sequence that $H^1\Big(\Ga, \Sel^{\s}(E/K_{\infty})(s)\Big)$ is a cofree $\Zp$-module. By \cite[Proposition 5.3.19(i)]{NSW}, this in turn implies that
    $\Big(\Sel^{\s}(E/K_{\infty})(s)\Big)^{\vee}$ has no nonzero finite $\La$-submodule. It then follows from this that $X^{\s}(E/K_{\infty})$ also has no nonzero finite $\La$-submodule.
\epf

\subsection{Non-primitive signed Selmer groups} \label{non-primitive Sel subsection}

For our purposes, it is useful to work with a non-primitive Selmer group (for instance, see \cite{G11}).
Let $\Si_0$ be a subset of $\Si_1$. The non-primitive signed Selmer group (with respect to $\Si_0$) is defined by
\[\Sel^{\s}_{\Si_0}(E/L_\infty) = \ker \Bigg(H^1(G_\Si(L_\infty),\Ep)\lra \bigoplus_{w\in \Si(L_\infty)-\Si_0(L_\infty)}J_w(E/L_{\infty}) \Bigg).\]
We shall write $X^{\s}_{\Si_0}(E/L_{\infty})$ for the Pontryagin dual of $\Sel^{\s}_{\Si_0}(E/L_\infty)$.

The non-primitive signed Selmer group and the signed Selmer group fit into the following commutative diagram
\[   \entrymodifiers={!! <0pt, .8ex>+} \SelectTips{eu}{}\xymatrix{
    0 \ar[r]^{} &\Sel^{\s}(E/L_{\infty})  \ar[d] \ar[r] &  H^1\big(G_{\Sigma}(L_{\infty}),\Ep\big)
    \ar@{=}[d] \ar[r] & \displaystyle \bigoplus_{u\in\Sigma(L_{\infty})}J_u(E/L_{\infty})\ar[d]^{}  \\
    0 \ar[r]^{} & \Sel^{\s}_{\Si_0}(E/L_{\infty}) \ar[r]^{} & H^1\big(G_{\Sigma}(L_\infty),\Ep\big) \ar[r] &
    \displaystyle\bigoplus_{u\in\Si(L_\infty)-\Si_0(L_\infty)}J_u(E/L_{\infty})  } \]
with exact rows. It then follows from a diagram chasing argument that we have the following exact sequence
\[0\lra \Sel^{\overrightarrow{s}}(E/L_{\infty}) \lra \Sel^{\overrightarrow{s}}_{\Si_0}(E/L_{\infty})\lra \bigoplus_{u\in \Si_0(F_{\infty})}H^1(L_{\infty,u}, \Ep). \]
We now record certain important properties of the non-primitive signed Selmer group.

\bp \label{no finite submodule}
Suppose that $\Sel^{\overrightarrow{s}}(E/K_{\infty})$ is a cotorsion $\Zp\ps{\Ga}$-module. Then the following statements are valid.

$(a)$ There are short exact sequences
\[0\lra \Sel^{\overrightarrow{s}}(E/K_{\infty}) \lra \Sel^{\overrightarrow{s}}_{\Si_0}(E/K_{\infty})\lra \bigoplus_{w\in \Si_0(K_{\infty})}H^1(K_{\infty,w}, \Ep)\lra 0 \] and
\[0\lra \Sel_{\Si_0}^{\overrightarrow{s}}(E/K_{\infty}) \lra H^1(G_\Si(K_\infty),\Ep)\lra \bigoplus_{w\in \Si(K_\infty)-\Si_0(K_\infty)}J_w(E/K_{\infty}) \lra 0. \]

$(b)$ $\theta\Big(X^{\s}(E/K_{\infty})\Big) = \theta\Big(X_{\Si_0}^{\s}(E/K_{\infty})\Big)$.

$(c)$ The cohomology groups
$H^i\Big(H, \Sel^{\s}_{\Si_0}(E/K_{\infty})\Big)$ are finite for every subgroup $H$ of $G$ and $i\geq 1$.

$(d)$ $X_{\Si_0}^{\s}(E/K_{\infty})$ has no non-trivial finite $\La$-submodule.
\ep

\bpf
 The first assertion is immediate from the above diagram and Proposition \ref{torsion surjective H2}. The second follows from a combination of assertion (a) and Lemma \ref{theta lemma2}, noting that $H^1(K_{\infty,u}, \Ep)$ is cofinitely generated over $\Zp$ by \cite[Proposition 2]{G89}. Assertion (c) can either be proven by a similar argument to that in Proposition \ref{tor cohomology}, or by taking $H$-invariant of the first short exact sequence in assertion (a) and applying the conclusion of Proposition \ref{tor cohomology}.

 Finally, we prove Assertion (d). From the first exact sequence in (a), we have
 \[ H^1\Big(\Ga, \Sel^{\overrightarrow{s}}(E/K_{\infty})\Big) \lra H^1\Big(\Ga, \Sel^{\overrightarrow{s}}_{\Si_0}(E/K_{\infty})\Big)\lra 0, \]
 where the zero on the right is a consequence of the well-known fact that
  \[H^1\left(\Ga,\bigoplus_{w\in \Si_0(K_{\infty})}H^1(K_{\infty,w}, \Ep) \right)=0.\]
  As seen in the proof of Proposition \ref{no finite submodules}, $H^1\Big(\Ga, \Sel^{\overrightarrow{s}}(E/K_{\infty})\Big)$ is cofree over $\Zp$. Hence so is $H^1\Big(\Ga, \Sel^{\overrightarrow{s}}_{\Si_0}(E/K_{\infty})\Big)$. It then follows from \cite[Proposition 5.3.19(i)]{NSW} that $X_{\Si_0}^{\s}(E/K_{\infty})$ has no non-trivial finite $\La$-submodule.
\epf

\subsection{Projectivity results} \label{main results}
Retain the notation of the previous subsections.
Let $v\in \Si_p^{ord}$, and $w$ a prime of $K$ lying above $v$. Write $k_w$ for the residue field of the local field $K_w$. Following Mazur \cite{Maz}, the prime $v$ is said to be
anomalous for $E/K$ if $|\widetilde{E}_v(k_w)|$ is divisible by $p$. Otherwise, we say that $v$ is non-anomalous
for $E/K$. Note that this definition is independent of the choice of $w$ above $v$.

As in \cite{G11}, we set
\[\Phi_{K/F} = \{ v\in \Si_1~|~\mbox{the inertia degree of $v$ in $K/F$ is divisible by $p$} \}. \]
We can now state and prove our results on the projectivity properties of the signed Selmer groups.

\bt \label{main theorem1}
Suppose that $(S1)-(S3)$ are satisfied and that $\Si_0$ contains $\Phi_{K/F}$. Assume that $X^{\s}(E/K_{\infty})$ is torsion over $\Zp\ps{\Ga}$ and that $\theta\Big(X^{\s}(E/K_{\infty})\Big)\leq 1$.

 Then $X^{\s}_{\Si_0}(E/K_{\infty})/X^{\s}_{\Si_0}(E/K_{\infty})[p]$ is quasi-projective as a $\Zp[G]$-module. Furthermore, if  $X^{\s}(E/K_{\infty})$ is finitely generated over $\Zp$, then
  $X^{\s}_{\Si_0}(E/K_{\infty})$ is quasi-projective as a $\Zp[G]$-module.
\et

\bpf
 Let $C=PQ$ be a cyclic subgroup of $G$ and $\e$ a $1$-dimensional character of $Q$. By Proposition \ref{no finite submodule}(a), we have a short exact sequence
\[0\lra \Sel_{\Si_0}^{\overrightarrow{s}}(E/K_{\infty})^\e \lra H^1(G_\Si(K_\infty),\Ep)^\e\lra \bigoplus_{w\in \Si(K_\infty)\setminus\Si_0(K_\infty)}J_w(E/K_{\infty})^\e \lra 0. \]
As seen in the proof of \cite[Proposition 3.2.1]{G11}, we have
\[h_p\Big(H^1(G_\Si(K_\infty),\Ep)^\e\Big) =1 \]
and
\[h_p\Big(J_w(E/K_{\infty})^\e\Big) =1 \]
for $w\in \Si^{ord}(K_\infty)\cup \big(\Si_1(K_\infty)\setminus\Si_0(K_\infty)\big)$. By Corollary \ref{H-cohomology of supersingular2} (and noting Remark \ref{relative +-}), the remaining local summands at the supersingular primes also have trivial $P$-Herbrand quotients. Hence we conclude that
\[h_p\Big(\Sel_{\Si_0}^{\overrightarrow{s}}(E/K_{\infty})^\e\Big) =1. \]
Now combining this with Proposition \ref{tor cohomology}, we obtain the first assertion of the theorem. The second assertion of the theorem will follow from an application of Proposition \ref{proj2}.
\epf

\bt \label{main theorem2}
Suppose that $(S1)-(S3)$ are satisfied and that $\Si_0$ contains $\Phi_{K/F}$. Assume that  $\Sel^{\s}(E/K_{\infty})$ is cotorsion over $\Zp\ps{\Ga}$, and that for each $v\in \Si_p^{ord}$, either $v$ is non-anomalous for $E/K$ or $v$ is tamely ramified in $K/F$.
 Then $X^{\s}_{\Si_0}(E/K_{\infty})$ has a free resolution of length 1 as a $\La[G]$-module.

 Moreover, if one further assumes that $X^{\s}(E/K_{\infty})$ is finitely generated over $\Zp$, then $X^{\s}_{\Si_0}(E/K_{\infty})$ is projective as a $\Zp[G]$-module.
 \et

\bpf
We first show that $H^i(H, \Sel^{\s}_{\Si_0}(E/K_{\infty}))=0$ for $i = 1,2$ and every subgroup $H$ of $G$. Let $L$ be the fixed field of $H$. Then we have the following commutative diagram
\[   \entrymodifiers={!! <0pt, .8ex>+} \SelectTips{eu}{}\xymatrix{
    0 \ar[r]^{} &\Sel^{\s}_{\Si_0}(E/L_{\infty})  \ar[d] \ar[r] &  H^1\big(G_{\Sigma}(L_{\infty}),\Ep\big)
    \ar[d]^{h} \ar[r] & \displaystyle \bigoplus_{v\in\Sigma(L_{\infty})\setminus\Si_0(L_\infty)}J_v(E/L_{\infty})\ar[d]^{g=\oplus g_v}  \ar[r] & 0\\
    0 \ar[r]^{} & \Sel^{\s}_{\Si_0}(E/K_{\infty})^H \ar[r]^{} & H^1\big(G_{\Sigma}(K_\infty),\Ep\big)^H \ar[r] &
    \left(\displaystyle\bigoplus_{w\in\Si(K_\infty)\setminus\Si_0(K_\infty)}J_w(E/K_{\infty})\right)^H&  } \]
with exact rows, where we note that the rightmost map on top sequence is surjective by Remark \ref{tor remark}. As seen in the proof of Proposition \ref{tor cohomology}, the map $h$ is an isomorphism and $H^i(H,H^2(G_{\Si}(K_{\infty}),\Ep))=0$ for $i\geq 1$. Hence one can conclude that
\[ H^i\Big(H,\Sel^{\s}_{\Si_0}(E/K_{\infty})\Big)\cong \begin{cases} \coker g,  & \mbox{if $i=1$}, \\
H^1\left(H,\displaystyle\bigoplus_{w\in\Si(K_\infty)\setminus\Si_0(K_\infty)}J_w(E/K_{\infty})\right), & \mbox{if $i= 2$}.\end{cases} \]
Note that $g=\oplus_v g_v$, where
\[ g_v:J_v(E/L_{\infty}) \lra \left(\bigoplus_{w|v} J_w(E/K_{\infty})\right)^H.\]

For good ordinary primes above $p$ and primes above $\Si_1\setminus\Si_0$, it was seen in the proof of \cite[Proposition 3.1.1]{G11} that $\coker g_v=0$ and \[H^1\left(H,~\bigoplus_{w|v}J_w(E/K_{\infty})\right)=0.\]
For the supersingular primes above $p$, the vanishing follows from Proposition \ref{H-cohomology of supersingular} and Remark \ref{relative +-}. In conclusion, one has $H^i(H, \Sel^{\s}_{\Si_0}(E/K_{\infty}))=0$ for $i=1,2$. Combining this observation with Propositions \ref{proj3} and \ref{no finite submodule}(d), we obtain the first assertion of the theorem.

Now supposing further that $X^{\s}(E/K_{\infty})$ is finitely generated over $\Zp$. Then, by Proposition \ref{no finite submodule}(a) and (d), we have that $X^{\s}_{\Si_0}(E/K_{\infty})$ is  finitely generated and free over $\Zp$. Combining this with the above established fact that $H^i(H, \Sel^{\s}_{\Si_0}(E/K_{\infty}))=0$ for $i=1$ and $i=2$, the second assertion then follows from Proposition \ref{proj1}.
\epf

\section{Applications of Theorems \ref{main theorem1} and \ref{main theorem2}} \label{application section}

In this section, we discuss some applications of Theorems \ref{main theorem1} and \ref{main theorem2}.

\subsection{Kida formula for signed Selmer groups} \label{Kida formula section}

As a start, we shall prove the Kida formula for our signed Selmer groups. Recall that
\[\Phi_{K/F} = \{ v\in \Si_1~|~\mbox{the inertia degree of $v$ in $K/F$ is divisible by $p$} \}. \]

\bl \label{p roots local}
Suppose that $K/F$ is a finite Galois $p$-extension.
For each $v\in \Phi_{K/F}$, we have $\mu_p\subseteq F_v$ and $\mu_p\subseteq K_w$ for every prime $w$ of $K$ above $v$.
\el

\bpf
 Let $w$ be a prime of $K$ above $v$, and $u$ a prime of $K_\infty$ above $w$ (and hence above $v$). Write $x$ for the prime of $F_\infty$ below $u$. By the assumption $v\in\Phi$, it follows that $K_w/F_v$ is ramified. As $v$ does not divide $p$, the field $F_{\infty,x}$ is a unramified $\Zp$-extension of $F_v$. Therefore, we have $K_w\nsubseteq F_{\infty,x}$. Recall that by \cite[Proposition 7.5.9]{NSW}, if a local field with residue characteristic $\neq p$ does not contain $\mu_p$, then its the maximal pro-$p$ extension is a $\Zp$-extension. But from the above discussion, we have seen that $K_{\infty, u}$ is an extension of $F_v$ with Galois group (slightly) bigger than $\Zp$. Hence we must have $\mu_p\subseteq F_v$.
  \epf

We can now prove the Kida's formula for the signed Selmer groups (compared with \cite{G11, HM, HSh, HL, PW}). Recall that we write $\la(M)$ for the $\la$-invariant of a torsion $\La$-module $M$ and that our conclusion is more in line with that in \cite{HSh}.

\bp \label{Kida formula}
Suppose that $p\geq 5$, that $(S1)-(S3)$ are satisfied and that $K/F$ is a finite Galois $p$-extension.
Furthermore, assume that $X^{\s}(E/K_{\infty})$ is torsion over $\La$ with $\theta\Big(X^{\s}(E/K_{\infty})\Big) \leq 1$. Then we have
\[\la\Big(\Sel^{\s}(E/K_{\infty})\Big) = |K:F|\cdot\la\Big(\Sel^{\s}(E/F_{\infty})\Big)+ \sum_{w\in P_1} (e_w-1) + 2\sum_{w\in P_2} (e_w-1),\]
where $P_1$ is the set of primes in $\Phi(K_{\infty})$ at which $E$ has split multiplicative reduction,
$P_2$ is the set of primes in $\Phi(K_{\infty})$  at which $E$ has good reduction with $E(K_{\infty,v})[p] \neq 0$, and $e_w$ is the ramification index of $w$ in $K_\infty/F_\infty$.
\ep

\bpf
  Set $\Si_0=\Phi_{K/F}$. By Theorem \ref{main theorem1}, the $\Zp[G]$-module $X^{\s}_{\Si_0}(E/K_{\infty})/X^{\s}_{\Si_0}(E/K_{\infty})[p]$ is  quasi-projective. It then follows from Corollary \ref{La proj cor} that
 \[\la\Big(\Sel^{\s}_{\Si_0}(E/K_{\infty})\Big) = |K:F|\cdot\la\Big(\Sel^{\s}_{\Si_0}(E/K_{\infty})^G\Big) =|K:F|\la\Big(\Sel^{\s}_{\Si_0}(E/F_{\infty})\Big), \]
 where the second equality follows from Lemma \ref{descent lemma}.
 Now, taking Proposition \ref{no finite submodule}(a) into account, we have
 \[\la\Big(\Sel^{\s}(E/K_{\infty})\Big) = |K:F|\cdot\la\Big(\Sel^{\s}(E/F_{\infty})\Big) + |K:F|\sum_{v\in \Phi_{K/F}(F_{\infty})} \corank_{\Zp}\Big(H^1(F_{\infty,v},\Ep\Big)\]
 \[  -
 \sum_{w\in \Phi_{K/F}(K_{\infty})} \corank_{\Zp}\Big(H^1(K_{\infty,w},\Ep\Big).\]
  Denote by $Q_1$ the set of primes in $\Phi_{K/F}(F_{\infty})$ at which $E$ has split multiplicative reduction, and
denote by $Q_2$ the set of primes in $\Phi_{K/F}(F_{\infty})$  at which $E$ has good reduction with $E(F_{\infty,v})[p] \neq 0$. Taking our assumption that $p\geq 5$ into account, it is straightforward to verify that $Q_i(K_{\infty})=P_i$ for $i=1,2$.
 Now, by a combination of \cite[Proposition 2]{G89} and \cite[Proposition 5.1]{HM}, we see that
 \[\corank_{\Zp}\Big(H^1(F_{\infty,v},\Ep\Big)=\begin{cases} 1,  & \mbox{if $v\in Q_1$}, \\
2, & \mbox{if $v\in Q_2$}, \\
0, & \mbox{if $v\in \Si_0(F_\infty)-Q_1\cup Q_2$},\end{cases} \]
and \[\corank_{\Zp}\Big(H^1(K_{\infty,w},\Ep\Big)=\begin{cases} 1,  & \mbox{if $w\in P_1$}, \\
2, & \mbox{if $w\in P_2$}, \\
0, & \mbox{if $w\in\Si_0(K_\infty)-P_1\cup P_2$},\end{cases} \]
where we note that to apply the results of Hachimori-Matsuno \cite[Proposition 5.1]{HM}, one requires the corresponding local fields to contain $\mu_p$ which is guaranteed by our Lemma \ref{p roots local}.
Therefore, we obtain
\[\la\Big(\Sel^{\s}(E/K_{\infty})\Big) = |K:F|\cdot\la\Big(\Sel^{\s}(E/F_{\infty})\Big) + |K:F|\sum_{v\in Q_1}1 + |K:F|\sum_{v\in Q_2}2 \]
 \[  -
 \sum_{w\in P_1}1 - \sum_{w\in P_2}2.\]
Since the set of primes of $K_{\infty}$ above $Q_i$ is precisely $P_i$ as noted above, we have
\[ |K:F|\sum_{v\in Q_i}1 = \sum_{v\in Q_i}\sum_{w|v}e_w = \sum_{w\in P_i}e_w\]
for $i=1,2$. Putting this back into the above equation, we obtain the equality of the theorem.
\epf

\br
Retaining the hypothesis of Theorem \ref{Kida formula} and assuming further that $\Sel^{\s}(E/F_{\infty})$ is cofinitely generated over $\Zp$, it then follows from Lemma \ref{descent lemma} that $\Sel^{\s}(E/K_{\infty})^G$ is cofinitely generated over $\Zp$. Since $G$ is a finite $p$-group, Nakayama lemma tells us that $\Sel^{\s}(E/K_{\infty})$ is cofinitely generated over $\Zp[G]$, and hence over $\Zp$. In this context, Theorem \ref{main theorem1} yields
\[\corank_{\Zp}\Big(\Sel^{\s}(E/K_{\infty})\Big) = |K:F|\corank_{\Zp}\Big(\Sel^{\s}(E/F_{\infty})\Big)+ \sum_{w\in P_1} (e_w-1) + 2\sum_{w\in P_2} (e_w-1),\]
which is the usual Kida's formula under the $\mu=0$ setting (see \cite{G11, HM, HL, PW}).
\er

\subsection{Characteristic elements of non-primitive signed Selmer groups}

We retain notation from Subsection \ref{char subsection} and Section \ref{Selmer}. Now if  $\Sel^{\s}(E/K_{\infty})$ is cotorsion over $\Zp\ps{\Ga}$, we can then attach characteristic elements to $X^{\s}(E/K_{\infty})$ and $X^{\s}_{\Si_0}(E/K_{\infty})$.
We shall write $\xi_{E}$ (resp., $\xi_{E,\Si_0}$) for a characteristic element of $X^{\s}(E/K_{\infty})$ (resp., $X^{\s}_{\Si_0}(E/K_{\infty})$). Recall that these characteristic elements live in the relative $K$-group $K_0(\La[G], Q_\La[G])$. It has been conjectured that these elements are related to a (conjectural) $p$-adic $L$-function in a precise manner (see \cite{CFKSV, Leip-adic, NP}). These $p$-adic $L$-functions are believed to satisfy certain integrality properties. Therefore, in view of these conjectures, one would expect that the characteristic elements satisfy certain integrality properties in an appropriate sense.
Motivated by this, our result is as follows (compare with \cite[Theorem 1]{NP}).

\bp \label{main theorem:char}
Suppose that $(S1)-(S3)$ are satisfied and that $\Si_0$ contains $\Phi$. Assume that $\Sel^{\s}(E/K_{\infty})$ is cotorsion over $\Zp\ps{\Ga}$ and that for each $v\in \Si_p^{ord}$, the prime $v$ is either non-anomalous for $E/K$ or tamely ramified in $K/F$.
Then $\xi_{E,\Si_0}$ belongs to the image of the following natural map
\[M_{\La[G]}\cap Q_{\La}[G]^\times \lra K_1( Q_\La[G]).\]
\ep

\bpf
This follows from a combination of Proposition \ref{char proj 1} and Theorem \ref{main theorem2}.
\epf

It is natural to ask whether we can obtain a similar conclusion for $\xi_{E}$. Of course, if $G$ is a group of order coprime to $p$, then we clearly have the conclusion. Therefore, we shall only concern ourselves with the situation when $p$ divides $|G|$. Unfortunately, at this point of writing, we do not have a complete answer to this. But we shall describe a situation, where one can prove such an integrality property. Recall from Theorem \ref{Kida formula}, $P_1$ is the set of primes in $\Phi_{K/F}(K_{\infty})$ at which $E$ has split multiplicative reduction, and $P_2$ is the set of primes $w$ in $\Phi_{K/F}(K_{\infty})$  at which $E$ has good reduction with $E(K_{\infty,w})[p] \neq 0$.

\bp \label{main theorem:char2}
Retain the settings of Theorem \ref{main theorem:char}. Suppose further that the sets $P_1$ and $P_2$ are empty.
Then $\xi_{E}$ belongs to the image of the following natural map
\[M_{\La[G]}\cap Q_{\La}[G]^\times \lra K_1( Q_\La[G]).\]
\ep

\bpf
Set $\Si_0=\Phi_{K/F}$.
By Proposition \ref{no finite submodule}(b), we have the following short exact sequence
\[0\lra \Sel^{\overrightarrow{s}}(E/K_{\infty}) \lra \Sel^{\overrightarrow{s}}_{\Si_0}(E/K_{\infty})\lra \bigoplus_{w\in \Si_0(K_{\infty})}H^1(K_{\infty,w}, \Ep)\lra 0. \]
Since the sets $P_1$ and $P_2$ are empty, it follows from the proof of Theorem \ref{Kida formula} that $H^1(K_{\infty,w}, \Ep)$ are finite for $w\in \Si_0(K_{\infty})$. Thus, the classes of $X^{\s}(E/K_{\infty})$ and $X^{\s}_{\Si_0}(E/K_{\infty})$ agree in $K_0(\La[G],Q_\La[G])$, and so they share the same characteristic elements. The conclusion of the proposition now follows from Theorem \ref{main theorem:char}.
\epf


\footnotesize
\begin{thebibliography}{00}

\bibitem{AL} S. Ahmed and M. F. Lim, On the Euler characteristics of signed Selmer groups, Bull. Aust. Math. Soc.  101 (2020), no. 2, 238-246.

\bibitem{AL2} S. Ahmed and M. F. Lim, On the signed Selmer groups of congruent elliptic curves with semistable reduction at all primes above $p$, Acta Arithmetica 197 (2021), no. 4, 353-377.

\bibitem{BL} K. B\"{u}y\"{u}kboduk and A. Lei, Integral Iwasawa theory of Galois representations
for non-ordinary primes. Math. Z. 286 (2017), no. 1-2, 361-398.

\bibitem{CFKSV} J.\ Coates, T.\ Fukaya, K.\ Kato, R.\ Sujatha and O.\ Venjakob, The $GL_2$ main conjecture for elliptic curves without complex multiplication. Publ. Math. Inst. Hautes \'Etudes Sci. No. 101, (2005), 163-208.

\bibitem{CG} J. Coates and R. Greenberg, Kummer theory for abelian varieties over
local fields. Invent. Math. 124 (1996), no. 1-3, 129-174.

\bibitem{EN} N. Ellerbrock and A. Nickel, On formal groups and
Tate cohomology in local fields. Acta Arith. 182 (2018), no. 3, 285-299.

\bibitem{G89} R. Greenberg, Iwasawa theory for $p$-adic representations, in: Algebraic number theory, 97-137, Adv. Stud. Pure Math., 17, Academic Press, Boston, MA, 1989.

\bibitem{G99} R. Greenberg, Iwasawa theory for elliptic curves, in: Arithmetic theory of elliptic curves (Cetraro, 1997), 51-144, Lecture Notes in Math., 1716, Springer, Berlin, 1999.

\bibitem{G11} R. Greenberg, Iwasawa theory, projective modules and modular representations. Mem. Amer. Math. Soc. 211 (2011), no. 992, vi+185 pp.


\bibitem{HM} Y. Hachimori and K. Matsuno,  An analogue of Kida's formula for the Selmer groups of elliptic curves. J. Algebraic Geom. 8 (1999), no. 3, 581--601.

\bibitem{HSh} Y. Hachimori and R. Sharifi, On the failure of pseudo-nullity of Iwasawa modules. J. Algebraic Geom. 14 (2005), no. 3, 567-591.


\bibitem{HL} J. Hatley and A. Lei, Arithmetic properties of signed Selmer groups at non-ordinary primes. Ann. Inst. Fourier (Grenoble) 69 (2019), no. 3, 1259-1294.

\bibitem{Iw81} K. Iwasawa, Riemann-Hurwitz formula and $p$-adic Galois representations for number fields.
Tohoku Math. J. (2) 33 (1981), no. 2, 263-288.

\bibitem{Kak11} M. Kakde, Proof of the main conjecture of noncommutative Iwasawa theory for totally
real number fields in certain cases. J. Algebraic Geom. 20 (2011), no. 4, 631-683.

\bibitem{Kak13} M. Kakde, The main conjecture of Iwasawa theory for totally
real fields. Invent. Math. 193 (2013), no. 3, 539-626.

\bibitem{K} K. Kato, $p$-adic Hodge theory and values of zeta functions of
modular forms, in: Cohomologies $p$-adiques et applications
arithm\'etiques. III., Ast\'erisque 295, 2004, ix, pp.
117-290.

\bibitem{Kim07}  B. D. Kim, The parity conjecture for elliptic curves at supersingular reduction primes. Compos. Math. 143 (2007), no. 1, 47-72.

\bibitem{KimPM} B. D. Kim, The plus/minus Selmer groups for supersingular primes. J. Aust. Math. Soc. 95 (2) (2013) 189-200.

\bibitem{Kim14} B. D. Kim, Signed-Selmer groups over the $\Zp^2$-extension of an imaginary quadratic
field. Canad. J. Math. 66 (2014), no. 4, 826-843.

\bibitem{KO} T. Kitajima and R. Otsuki,
On the plus and the minus Selmer groups
for elliptic curves at supersingular primes. Tokyo J. Math. 41 (2018), no. 1, 273-303.

\bibitem{Kob} S. Kobayashi, Iwasawa theory for elliptic curves at supersingular primes. Invent. Math. 152 (2003), no. 1, 1-36.

\bibitem{Leip-adic} A. Lei, Non-commutative $p$-adic $L$-functions for supersingular primes. Int. J. Number Theory 8 (2012), no. 8, 1813--1830

\bibitem{LLMordell}  A. Lei and M. F. Lim, Mordell-Weil ranks and Tate-Shafarevich groups of elliptic curves with mixed-reduction type over cyclotomic extensions, arXiv:1911.10643 [math.NT].

\bibitem{LLAkashi} A. Lei and M. F. Lim,  Akashi series and Euler characteristics of signed Selmer groups of elliptic curves with semistable reduction at primes above
$p$, accepted for publication in J. Theor. Nombres Bordeaux.

\bibitem{LZ} A. Lei and S. Zerbes, Signed Selmer groups over $p$-adic Lie extensions. J. Th\'eor. Nombres
Bordeaux 24 (2012), no. 2, 377-403.

\bibitem{Maz} B. Mazur, Rational points of abelian varieties with values in
towers of number fields. Invent. Math. 18 (1972), 183-266.

\bibitem{NSW} J. Neukirch, A. Schmidt and K. Wingberg,
Cohomology of Number Fields, 2nd edn., Grundlehren Math.
Wiss. 323 (Springer-Verlag, Berlin, 2008).

\bibitem{NP} A. Nichifor and B. Palvannan, On free resolutions of Iwasawa modules. Doc. Math. 24 (2019), 609-662.

\bibitem{No} E. Noether, Normalbasis bei K\"{o}rpern ohne h\"ohere Verzweigung. J. Reine Angew. Math. 167 (1932), 147-152.

\bibitem{PW} R. Pollack and T. Weston, Kida's formula and congruences. Doc. Math. 2006, Extra Vol., 615-630.

\bibitem{Re} I. Reiner, Maximal orders,  Corrected reprint of the 1975 original. With a foreword by M. J. Taylor. London Mathematical Society Monographs. New Series, 28. The Clarendon Press, Oxford University Press, Oxford, 2003. xiv+395 pp.

\bibitem{Wi} M. Witte, On a localisation sequence for the $K$-theory of skew power series rings. J. K-Theory 11 (2013), no. 1, 125-154.
\end{thebibliography}
\end{document}